\documentclass[AMA,STIX1COL,doublespace]{WileyNJD-v2}
\usepackage{array}
\usepackage{stmaryrd}
\usepackage[caption=false]{subfig}
\usepackage[flushleft]{threeparttable}
\usepackage{multirow}
\allowdisplaybreaks

\articletype{Draft}%

\begin{document}

\title{A hybridizable discontinuous Galerkin method for simulation of electrostatic problems with floating potential conductors}

\author[1]{Liang Chen*}

\author[1]{Ming Dong}

\author[2]{Ping Li}

\author[1]{Hakan Bagci}

%\authormark{AUTHOR ONE \textsc{et al}}

\address[1]{\orgdiv{Division of Computer, Electrical, and Mathematical Science and Engineering (CEMSE)}, \orgname{King Abdullah University of Science and Technology (KAUST)}, \orgaddress{\state{Thuwal 23955-6900}, \country{Saudi Arabia}}}

\address[2]{\orgdiv{Department of Electrical Engineering}, \orgname{Shanghai Jiao Tong University}, \orgaddress{\state{Shanghai, 200240}, \country{China}}}

\corres{*Liang Chen. \email{liang.chen@kaust.edu.sa}}

%\presentaddress{This is sample for present address text this is sample for present address text}

\abstract[Summary]{
In an electrostatic simulation, an equipotential condition with an undefined/floating potential value has to be enforced on the surface of an isolated conductor. If this conductor is charged, a nonzero charge condition is also required. While implementation of these conditions using a traditional finite element method (FEM) is not straightforward, they can be easily discretized and incorporated within a discontinuous Galerkin (DG) method. However, DG discretization results in a larger number of unknowns as compared to FEM. In this work, a hybridizable DG (HDG) method is proposed to alleviate this problem. Floating potential boundary conditions, possibly with different charge values, are introduced on surfaces of each isolated conductor and are weakly enforced in the global problem of HDG. The unknowns of the global HDG problem are those only associated with the nodes on the mesh skeleton and their number is much smaller than the total number of unknowns required by DG. Numerical examples show that the proposed method is as accurate as DG while it improves the computational efficiency significantly.}


%\jnlcitation{\cname{%
%		\author{L. Chen}, 
%		\author{M. Dong},
%		\author{P. Li}, and
%		\author{H. Bagci} } (\cyear{2020}), 
%\ctitle{A hybridizable discontinuous Galerkin method for simulation of electrostatic problems with floating potential conductors}, \cjournal{INT. J. NUMER. METH. ENG.}, \cvol{2020;00:1--8}.}

\maketitle

%\footnotetext{\textbf{Abbreviations:} ANA, anti-nuclear antibodies; APC, antigen-presenting cells; IRF, interferon regulatory factor}

\section{Introduction}\label{sec1}
Isolated conductors exist in a wide range of electrical and electronic systems, such as electrode cores of high-voltage inductors~\cite{Aracil2014}, metallic separators of IEC surge arresters~\cite{Andjelic2018}, defects in ultra-high-voltage gas-insulated switchgear~\cite{Zeng2018}, passive electrodes of earthing systems~\cite{Zildzo2009}, conductors of floating-gate transistors~\cite{Kahng1967}, and, more recently, metallic nanostructures extensively used in optoelectronic devices~\cite{Park2012}. In electrostatic simulations of these systems, these conductors result in equipotential surfaces with unfixed (i.e., floating) electric potential values (which depend on the simulation parameters and the geometries of the structures involved) and are referred to as floating potential conductors (FPCs).

Even though execution of these electrostatic simulations by using a finite element method (FEM), which solves the Poisson equation, has become a common practice, accurate and efficient incorporation of FPC models within a FEM framework is still not a trivial task. In recent years, various techniques have been introduced to address this challenge ~\cite{Konrad1996, Blaszczyk1994, Takuma1997, Dong2010, Fu2012, Rincon2018, Dong2018, Aiello2016}. The most commonly used methods among these are the virtual permittivity method (VPM)~\cite{Konrad1996}, the matrix reduction method (MRM)~\cite{Dong2010}, and the charge simulation method (CSM)~\cite{Blaszczyk1994, Takuma1997, Dong2018}. Each of methods has pros and cons regarding the accuracy and efficiency of the solution, the ability to account for charges on FPCs, and the ease of implementation~\cite{Fu2012, Rincon2018, Chen2020float}. VPM uses a dielectric material with a very high ``virtual'' permittivity to approximate the conductor. This method is very straightforward to implement since it does not require any modifications to be done on an existing FEM code. But its accuracy depends on the value of the virtual permittivity. Accurate representation of a conductor requires a very high virtual permittivity value but this, in return, makes the FEM matrix ill-conditioned and leads to a less accurate solution. MRM produces more accurate results but it requires rather significant modifications to be done on the original FEM code~\cite{Dong2010, Fu2012, Rincon2018}. In addition, when an FPC is charged, a nonzero charge condition must be imposed on its surface. Both VPM and MRM can not account for this nonzero charge condition~\cite{Dong2010, Fu2012, Rincon2018}. CSM can account for charge conditions since it enforces a specific charge distribution on an FPC but this requires a priori knowledge of simulation results or multiple iterative simulations~\cite{Dong2010, Fu2012, Rincon2018, Dong2018}. 
Several boundary element methods (BEMs) have also been developed for modeling FPC in electrostatic simulations~\cite{Andjelic2018, Amann2014}. In~\cite{Andjelic2018} a total electric charge condition is applied to determine the potential of uncharged FPCs. BEM is often preferred over FEM for unbounded problems with homogeneous or piece-wise homogeneous materials. In~\cite{Amann2014}, the Poincare-Steklov operator is used to enforce constraints corresponding to the floating potential.

Recently, it has been shown that FPCs can easily be accounted for using the discontinuous Galerkin (DG) method~\cite{Chen2020float}. By weakly imposing a so-called floating potential boundary condition (FPBC) through the numerical flux, DG can accurately model FPCs with non-zero charge conditions. The implementation of this approach in existing DG codes is rather straightforward. In addition, by enforcing FBPC on the surfaces of FPCs, the requirement to discretize their volumes is removed, which reduces the total number of unknowns. However, even with this reduction, the number of unknowns required by DG is still larger than that of the traditional FEM, which might lead to a considerable increase in computational cost depending on the problem being analyzed. 

In recent years, the hybridizable DG (HDG)~\cite{Cockburn2009unified} method has drawn a lot of attention. HDG addresses the fundamental weakness of DG, i.e., reduces the total number of knowns by applying a static condensation technique within the DG framework~\cite{Cockburn2016static, Cockburn2016hdg}. In HDG, a hybrid variable is introduced on the mesh skeleton and mesh elements exchange information only through this variable. This approach allows for locating globally coupled degrees of freedom only on the mesh skeleton and results in a global matrix system (in the unknown hybrid variable) with a dimension much smaller than that of the matrix system generated by DG. The local unknowns are then recovered using the hybrid variable that is obtained by solving this smaller global system. Furthermore, HDG can achieve superconvergence by applying a local post-processing technique where the solution represented using polynomials of order $p$ leads to a convergence of order $p+2$ in accuracy. HDG has been applied to various problems~\cite{Cockburn2016static, Cockburn2016hdg}, such as fluid dynamics~\cite{Nguyen2009implicit, Nguyen2009implicitnonlinear}, and acoustic, elastic, and electromagnetic wave propagations~\cite{Nguyen2011high, Nguyen2011hybridizable, Li2013hybridizable}. Comparative studies~\cite{Kirby2012ToCGorHDG2D, Yakovlev2016ToCGorHDG3D} between HDG and FEM show that HDG outperforms FEM in efficiency for some cases, while it retains other advantages of DG~\cite{Kirby2012ToCGorHDG2D, Yakovlev2016ToCGorHDG3D, Cockburn2016static}.

In this work, an HDG-based framework is proposed to implement FPBC in electrostatic simulations. The local problem is formulated with a Dirichlet boundary condition and the electric potential is chosen as the hybrid variable in the global problem. The floating potential values on each FPC are also left as unknowns in the global problem. The dimension of the resulting global system is equal to the number of nodes on the mesh skeleton (where the hybrid variable is defined) plus the number of  FPCs, which in total is significantly smaller than the dimension of the matrix system that would be solved by DG. Other advantages of using the FPBC, such as the ease of implementation, the high solution accuracy, and the ability to account for non-zero charge conditions are inherited by this HDG-based framework. Table~\ref{comparison} compares the properties of this proposed method to other methods briefly described above. Note that the low efficiency of DG is because of the larger number of unknowns, while for CSM it is due to multiple-simulation requirement~\cite{Rincon2018}.

The rest of the paper is organized as follows. Section II starts with the mathematical model of the electrostatic problem in the presence of FPCs, then it presents the HDG formulation for this problem, including the hybridized strong form, the weak form, and the discretized matrix system. Section III presents several numerical examples and Section IV provides a summary.

\begin{table}[!b]
	\renewcommand{\arraystretch}{1.05}
	\centering
	\caption{Comparision of FEM-type methods used for modeling FPCs in electrostatic simulations.}
	\label{comparison}
	\setlength{\tabcolsep}{3pt}
	\begin{tabular}{p{200pt} p{44pt} p{44pt} p{44pt} p{44pt} p{44pt}}
%	\begin{tabular}{c c c c c c}
		\hline
		& VPM & MRM & CSM  & DG & HDG\\ \hline
		Capability of modeling charged FPCs & no & no & yes & yes & yes \\ \hline
		Accuracy of representing physical conditions & low & high & low & high & high\\ \hline
		Easy adaption into existing code & easy & difficult & difficult & easy & easy\\ \hline
		Needing of internal meshes & yes & yes & yes & no & no \\ \hline
		Relative efficiency & moderate & high & low & low & moderate \\ \hline
	\end{tabular}
	\label{tab1}
\end{table}

\section{Formulation}\label{sec2}
\subsection{Mathematical Model}
\begin{figure}[!t]
	\centerline
	{\includegraphics[width=0.4\columnwidth]{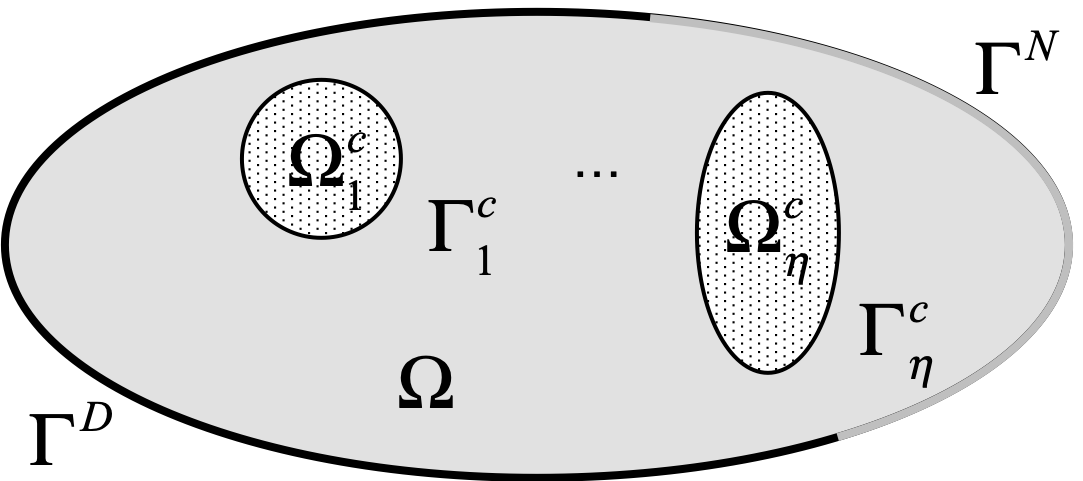}}
	\caption{Schematic description of an electrostatic problem involving multiple isolated conductors.}
	\label{Domain}
\end{figure}

Consider the electrostatic problem described in Figure~\ref{Domain}. $M$ isolated conductors $\Omega _1^c,\Omega _2^c,...,\Omega _M^c$ are distributed inside domain $\Omega $. Denote the surface and charge of each conductor by $\Gamma _\eta ^c$ and  $Q _\eta ^c$, $\eta  = 1,2,...M$, respectively and denote the domain boundary by $\partial\Omega$, $\partial\Omega=\Gamma^D\cup\Gamma^N\cup\Gamma^c$, where $\Gamma^c = \cup_{\eta=1}^M \Gamma _\eta ^c$, and $\Gamma ^D$ and $\Gamma ^N$ represent the boundaries where Dirichlet and Neumann boundary conditions are enforced, respectively. The electrostatic problem is described by the following boundary value problem (BVP)~\cite{Chen2020float}
\begin{alignat}{2}
\label{BVP0} & \nabla  \cdot [\varepsilon ({\mathbf{r}})\nabla \varphi ({\mathbf{r}})] =  - \rho ({\mathbf{r}}),
\quad && {\mathbf{r}} \in \Omega\\
\label{BVP1} & \varphi ({\mathbf{r}}) = {f^D}({\mathbf{r}}),
\quad && {\mathbf{r}} \in {\Gamma^D}\\
\label{BVP2} & \hat{\mathbf{n}}(\mathbf{r}) \cdot [\varepsilon ({\mathbf{r}})\nabla \varphi ({\mathbf{r}})] = {f^N}({\mathbf{r}}),
\quad && {\mathbf{r}} \in {\Gamma^N}\\
\label{BVP3} & \varphi ({\mathbf{r}}) = \varphi _\eta ^c, \quad \oint_{\Gamma _\eta ^c} {\hat{\mathbf{n}}(\mathbf{r}) \cdot [\varepsilon ({\mathbf{r}})\nabla \varphi ({\mathbf{r}})] d{\mathbf{r}}} = - Q _\eta ^c
\quad && {\mathbf{r}} \in \Gamma _\eta ^c.
\end{alignat}
In~\eqref{BVP0}-\eqref{BVP3}, $\varphi({\mathbf{r}})$ is the electric potential distribution to be solved for, $\varepsilon ({\mathbf{r}})$ is the permittivity, $\rho ({\mathbf{r}})$ is the charge density, ${f^D}({\mathbf{r}})$ and ${f^N}({\mathbf{r}})$ are the coefficients associated with the Dirichlet and Neumann boundary conditions, respectively, and $\hat{\mathbf{n}}(\mathbf{r})$ denotes the outward-pointing normal vector of the corresponding surface. Equation~$\eqref{BVP3}$ represents the physical conditions on FPCs. On each FPC, the equipotential value $\varphi _\eta ^c$ is an unknown, and the total charge $Q _\eta ^c$ is assumed known. The charge condition in~\eqref{BVP3}, i.e., the total electric flux is equal to the total charge, provides the constraint that makes $\varphi _\eta ^c$ unique.

\subsection{The Hybridizable Discontinuous Galerkin Method}
\subsubsection{The strong form}
To develop the HDG method, \eqref{BVP0}-\eqref{BVP3} are expressed as a first-order partial differential equation system by using the electric field ${\mathbf{E}}({\mathbf{r}}) =  - \nabla \varphi ({\mathbf{r}})$. The BVP becomes finding $\varphi({\mathbf{r}})$ and $\mathbf{E}({\mathbf{r}})$ such that
\begin{alignat}{2}
\label{BVP_0}  & \nabla  \cdot [\varepsilon ({\mathbf{r}}){\mathbf{E}}({\mathbf{r}})] = \rho ({\mathbf{r}}), \quad && {\mathbf{r}} \in \Omega\\
\label{BVP_1}  & {\mathbf{E}}({\mathbf{r}}) + \nabla \varphi ({\mathbf{r}}) = 0, \quad &&  {\mathbf{r}} \in \Omega\\
\label{BVP_2}  & \varphi ({\mathbf{r}}) = {f^D}({\mathbf{r}}), \quad && {\mathbf{r}} \in {\Gamma^D}\\
\label{BVP_3}  & - \hat{\mathbf{n}}(\mathbf{r}) \cdot [\varepsilon ({\mathbf{r}}){\mathbf{E}}({\mathbf{r}})] = {f^N}({\mathbf{r}}), \quad && {\mathbf{r}} \in {\Gamma^N}\\
\label{BVP_4}  & \varphi ({\mathbf{r}}) = \varphi _\eta ^c, \quad \oint_{\Gamma _\eta ^c} {{\hat{\mathbf{n}}}({\mathbf{r}}) \cdot [\varepsilon ({\mathbf{r}}){\mathbf{E}}({\mathbf{r}})] d{\mathbf{r}}}  = Q _\eta ^c,\quad && {\mathbf{r}} \in \;\Gamma _\eta ^c.
\end{alignat}

\begin{figure}[t]
	\centerline{\includegraphics[width=0.5\columnwidth]{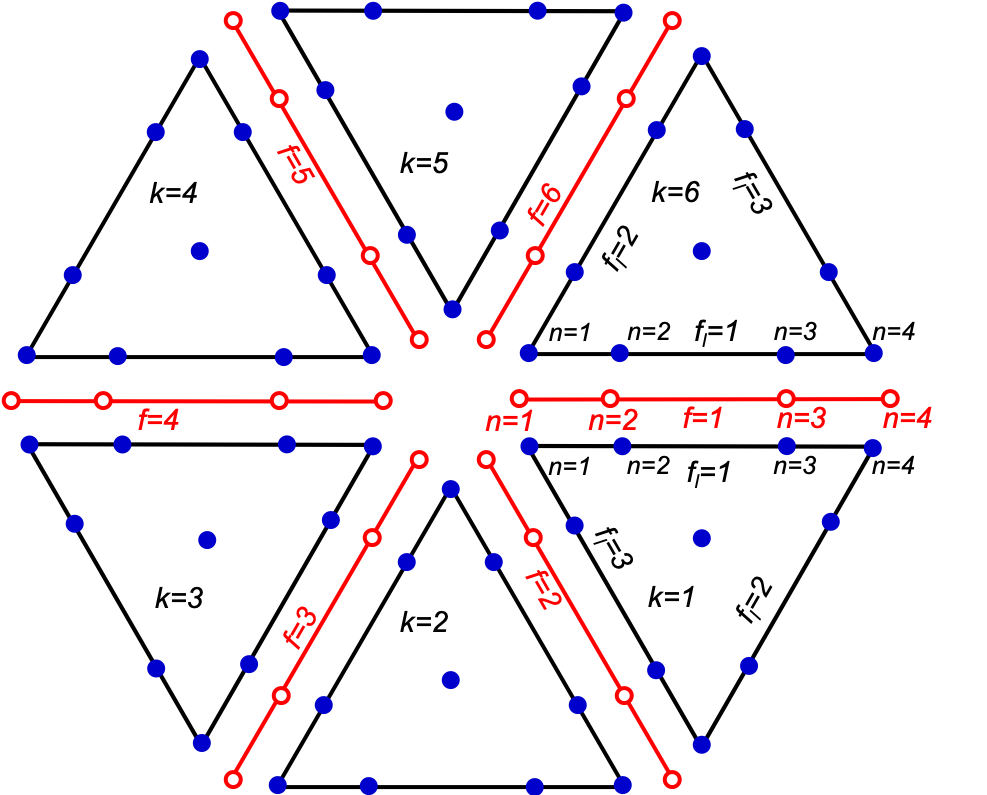}}
	\caption{Illustration of the local unknowns (blue dots) and the global unknowns (red circles) for HDG with $p=3$.\label{mesh}}
\end{figure}

Partition $\Omega$ into $K$ non-overlapping tetrahedrons, $\Omega = \cup _{k=1}^K {\Omega _k}$, and denote the surface of ${\Omega _k}$ and its the outward-pointing unit vector normal by $\partial {\Omega _k}$ and $\hat{\mathbf{n}}_k({\mathbf r})$, respectively. Further, denote the interior skeleton (see the red lines illustrated in Fig.~\ref{mesh})
% (edges and faces in two dimension (2D) and three dimension (3D), respectively) 
of this mesh of tetrahedrons by $\Gamma=[\cup_{k=1}^{K}\partial\Omega_k]\setminus\partial\Omega$ and denote the total number of faces on $\Gamma$ by $N_f$. The number of faces per each element is denoted by $N_{fe}$. Hereinafter, the explicit dependency on $\mathbf{r}$ is omitted for the brevity of the notation. Following the hybridization approach developed in~\cite{Cockburn2009unified}, the local problem on element $k$ is defined as
\begin{alignat}{2}
\label{Localfp0} & \nabla \cdot (\varepsilon {\mathbf{E}})_k = \rho,
\quad && \text{in } \Omega_k\\
\label{Localfp1} & {\mathbf{E}_k} + \nabla \varphi_k = 0,
\quad && \text{in } \Omega_k\\
\label{Localfp2} & \varphi_k = \hat{\varphi}_k,
\quad && \text{on } \partial\Omega_k \cap \Gamma \\
\label{Localfp5} & \varphi_k = \hat{\varphi}^c,
\quad && \text{on } \partial \Omega_k \cap \Gamma^c\\
\label{Localfp3} & \varphi_k = f^D,
\quad && \text{on } \partial\Omega_k \cap \Gamma^D \\
\label{Localfp4} & {\hat{\mathbf{n}}}_k \cdot (\varepsilon{\mathbf{E}})_k=f^N,
\quad && \text{on } \partial\Omega_k \cap \Gamma^N
\end{alignat}
where $\varphi_k$, $\mathbf{E}_k$, and $\hat{\varphi}_k$ are the local variables defined on element $k$, $(\varepsilon{\mathbf{E}})_k=\varepsilon_k {\mathbf{E}}_k$, and $\varepsilon_k$ is the permittivity, which is assumed constant in element $k$. Here, $\hat{\varphi}_k$ and $\hat{\varphi}^c$ are hybrid variables, which satisfy a global problem (to be described below). In particular, $\hat{\varphi}^c$ is the value of the floating potential on the boundary $\Gamma^c$. Without loss of generality, multiple  $\hat{\varphi}^c$ variables can be defined for multiple FPCs independently. Equations \eqref{Localfp0}-\eqref{Localfp4} describe a local BVP on each element. Once $\hat{\varphi}_k$ and $\hat{\varphi}^c$ are known, they can be used as Dirichlet boundary values to solve this local BVP. Note that this also means that the local variables $\varphi_k$ and $\mathbf{E}_k$ can be expressed as functions of the hybrid variables $\hat{\varphi}_k$ and $\hat{\varphi}^c$.
%The weak form of \eqref{Localfp0}-\eqref{Localfp4} 

The hybrid variables are required to satisfy the global problem defined using the transmission condition~\cite{Cockburn2009unified, Cockburn2016static}
\begin{align}
\label{Globalfp0} \llbracket {\hat{\mathbf{n}}_k} \cdot (\varepsilon{\mathbf{E}})_k \rrbracket & =0, \quad \text{on } \partial\Omega_k \cap \Gamma
\end{align}
and the charge condition~\cite{Chen2020steadystate}
\begin{align}
\label{Charge} \sum_{f=1}^{N_f^c} \hat{\mathbf{n}}_k \cdot (\varepsilon{\mathbf{E}})_k & = Q^c, \quad \text{on } \partial\Omega_k \cap \Gamma^c.
\end{align}
Here, $\llbracket \odot \rrbracket = \odot^+ + \odot^-$ defines the jump at the inter-element boundaries and $N_f^c$ is the total number of faces on the FPC. Since $\varphi_k$ and $\mathbf{E}_k$ are functions of the hybrid variables, \eqref{Globalfp0}-\eqref{Charge} can be cast into a global matrix system with unknowns $\hat{\varphi}_k$ and $\hat{\varphi}^c$.

The hybridized system described above is amenable to the static condensation of continuous finite element methods (CFEM) and the hybridization of mixed finite element methods (MFEM)~\cite{Cockburn2009unified, Cockburn2016static}. Different from CFEM and MFEM, HDG solves the local BVP by DG and enforces the transmission boundary condition weakly using the numerical flux~\cite{Cockburn2009unified, Cockburn2016static}. This gives rise to the generalization of the FPBC from DG, where the FPBC is weakly imposed using the numerical flux~\cite{Chen2020float}, to HDG.

\subsubsection{The weak form}
Let $\mathbb{P}^p$ denote the space of polynomial functions of degree at most $p$ $(p\geq 1)$, then the following discrete finite element spaces can be introduced
\begin{alignat*}{3}
\mathbb{L} & = \{ l \in L_2(\Omega); {}&&{} l|_{\Omega_k} \in \mathbb{P}^p(\Omega_k), {}&&{} \forall \Omega_k \in \Omega \} \\
\mathbb{W} & = \{ \mathbf{w} \in [L_2(\Omega)]^3;\quad{}&&{} \mathbf{w}|_{\Omega_k} \in [\mathbb{P}^p(\Omega_k)]^3,\quad {}&&{} \forall \Omega_k \in \Omega \} \\
\mathbb{V} & = \{ v \in L_2(\Gamma); {}&&{} v|_{\Gamma_f} \in \mathbb{P}^p(\Gamma_f), {}&&{} \forall \Gamma_f \in \Gamma\cup \partial \Omega \}.
\end{alignat*}
%\textbf{Do we need to define the interpolation of $\hat{\varphi}$ on $\Gamma$?}
Let $\left(\cdot,\cdot\right)_{\Omega}$ denote the $L_2$ inner product in the domain $\Omega$
\begin{equation*}
\left(u,v\right)_{\Omega} = \int_{\Omega}uv d\mathbf{r}, \quad 
\left(\mathbf{u},\mathbf{v}\right)_{\Omega} = \int_{\Omega}\mathbf{u}\mathbf{v} d\mathbf{r}
\end{equation*}
and $\langle \cdot,\cdot\rangle _{\Gamma}$ denote the $L_2$ inner product on the face $\Gamma$
\begin{equation*}
\langle u,v\rangle_{\Gamma} = \int_{\Gamma}uv d\mathbf{r}, \quad \langle\mathbf{u},\mathbf{v}\rangle_{\Gamma} = \int_{\Gamma}\mathbf{u}\mathbf{v} d\mathbf{r}.
\end{equation*}

Following the classical DG approach~\cite{Hesthaven2008, Cockburn2009unified, Chen2020float, Chen2020steadystate, Chen2019multiphysics}, the weak form of the local problem~\eqref{Localfp0}-\eqref{Localfp4} is defined as
\begin{align}
\label{PSweak0} - \left(  \nabla {l}, {(\varepsilon {\mathbf{E}})_k} \right)_{\Omega_k}  +
\langle {l, \hat{\mathbf{n}}_k \cdot {{(\varepsilon {{\mathbf{E}}})}_k^*}} \rangle_{\partial {\Omega _k}} & = \left( {l, \rho} \right)_{\Omega_k} \\
\label{PSweak1}  \left( {\mathbf{w}, \mathbf{E}_k}\right)_{\Omega_k} - \left( {\nabla \cdot \mathbf{w}, {\varphi _k} }\right)_{\Omega_k} + \langle { \hat{\mathbf{n}}_k \cdot \mathbf{w}, \varphi _k^* }\rangle _{\partial {\Omega_k}} & = 0
\end{align}
where $\varphi _k \in \mathbb{L}$, $\mathbf{E}_k \in \mathbb{W}$, and $\hat{\varphi}_f \in \mathbb{V}$ are the approximate solutions sought for (for the sake of simplicity, the same notations are used for the variables in the weak form and the strong form), and the weak form of the global problem~\eqref{Globalfp0}-\eqref{Charge} is defined as
\begin{align}
\label{PSweak2} \sum\limits_{k=1}^{K}{\langle v,\hat{\mathbf{n}}_k \cdot {({\varepsilon \mathbf{E}})_k^*}{{\rangle}_{\partial {{\Omega }_{k}} \cap \Gamma }}}& =0 \\
\label{PSweak3} \sum\limits_{k=1}^{K}{\langle v,\hat{\mathbf{n}}_k \cdot {({\varepsilon \mathbf{E}})_k^{*}} \rangle _{\partial {{\Omega }_{k}}\cap {\Gamma^c}}}& =Q^c
\end{align}
where the numerical fluxes $\varphi_k^*$ and $({\varepsilon \mathbf{E}})_k^*$ are chosen as
\begin{align}
\label{flux}
\left.
\begin{aligned}
& \varphi_k^* = \hat{\varphi}_k, \hspace{0.1cm}
&& {\hat{\mathbf{n}}}_k \cdot ({\varepsilon \mathbf{E}})_k^* = {\hat{\mathbf{n}}}_k \cdot ({\varepsilon \mathbf{E}})_k + \tau(\varphi_k-\hat{\varphi}_k), \hspace{0.1cm}
&& {\text{on }} \partial \Omega_k \cap \Gamma\\
& \varphi_k^* = f_k^D, \hspace{0.1cm}
&& {\hat{\mathbf{n}}}_k \cdot ({\varepsilon \mathbf{E}})_k^* = {\hat{\mathbf{n}}}_k \cdot ({\varepsilon \mathbf{E}})_k + \tau(\varphi_k-f_k^D), \hspace{0.1cm}
&& {\text{on }} \partial \Omega_k \cap \Gamma^D\\
& \varphi_k^* = \hat{\varphi}_k^c, \hspace{0.1cm}
&& {\hat{\mathbf{n}}}_k \cdot ({\varepsilon \mathbf{E}})_k^* = {\hat{\mathbf{n}}}_k \cdot ({\varepsilon \mathbf{E}})_k, \hspace{0.1cm}
&& {\text{on }} \partial \Omega_k \cap \Gamma^c\\
& \varphi_k^* = \varphi_k, \hspace{0.1cm}
&& {\hat{\mathbf{n}}}_k \cdot ({\varepsilon \mathbf{E}})_k^* = f_k^N, \hspace{0.1cm}
&& {\text{on }} \partial \Omega_k \cap \Gamma^N.
\end{aligned}
\right.
\end{align}
The numerical fluxes in~\eqref{flux} follow from the local DG (LDG) method~\cite{Cockburn1998} and the resulting HDG method is also called LDG-H method~\cite{Cockburn2009unified,Cockburn2016static}. The stabilization parameter $\tau$ is of order $1/h$~\cite{Cockburn2009unified,Cockburn2016static}, where $h$ is the element edge length.

Substituting \eqref{flux} into ~\eqref{PSweak0}-\eqref{PSweak3} gives
\begin{align}
\label{weakfp0} & 
\langle l,\tau \varphi_k \rangle_{\partial \Omega_k \setminus \Gamma^N} +
\left( l, \nabla \cdot (\varepsilon {\mathbf{E}})_k \right)_{\Omega_k} - 
\langle l,\hat{\mathbf{n}}_k\cdot (\varepsilon {\mathbf{E}})_k \rangle_{\partial \Omega_k \cap \Gamma^N} - 
\langle l,\tau \hat{\varphi}_k \rangle_{\partial \Omega_k \cap \Gamma} = 
\left(l,f_k\right)_{\Omega_k} + 
\langle l,\tau f_k^D \rangle_{\partial \Omega_k \cap \Gamma^D} -
\langle l,f_k^N \rangle _{\partial \Omega_k \cap \Gamma^N} \\
\label{weakfp1} & 
\left(\nabla\cdot\mathbf{w},\varphi_k\right)_{\Omega_k} -
\langle \hat{\mathbf{n}}_k\cdot \mathbf{w},\varphi_k \rangle_{\partial \Omega_k \cap \Gamma^N} -
\left( \mathbf{w}, (\varepsilon {\mathbf{E}})_k \right)_{\Omega_k} -
\langle \hat{\mathbf{n}}_k \cdot \mathbf{w},\hat{\varphi}_k \rangle_{\partial \Omega_k \cap \Gamma } -
{\langle {\hat{\mathbf{n}}_k} \cdot {\mathbf{w}},\hat{\varphi}_k^c \rangle _{\partial {\Omega _k} \cap {\Gamma^c}}} = 
\langle \hat{\mathbf{n}}_k \cdot \mathbf{w},f_k^D \rangle_{\partial \Omega_k \cap \Gamma^D} \\
\label{weakfp2} & 
\sum\limits_{k=1}^{K}{ \left\{
	\langle v,\tau \varphi_k \rangle _{\partial {{\Omega }_{k}}\cap \Gamma} + \langle v,\hat{\mathbf{n}}_k \cdot (\varepsilon {\mathbf{E}})_k \rangle _{\partial {{\Omega }_{k}}\cap \Gamma}-\langle v,\tau \hat{\varphi}_k \rangle _{\partial {{\Omega }_{k}}\cap \Gamma}
	\right\} }=0 \\
\label{weakfp3} & 
\sum\limits_{k=1}^{K}
{ \langle v,\hat{\mathbf{n}}_k \cdot (\varepsilon {\mathbf{E}})_k \rangle _{\partial {{\Omega }_{k}}\cap {\Gamma^c}} }=Q^c.
\end{align}
In \eqref{weakfp0}, the divergence theorem is applied once to symmetrize the problem. Equations \eqref{weakfp0}-\eqref{weakfp3} represents the final system to be solved for $\mathbf{E}_k$, $\varphi_k$ and $\hat{\varphi}_k$.

\subsubsection{The discrete system}
Using Lagrange polynomials~\cite{Hesthaven2008}, the nodal interpolation of $\varphi_k$ and $\mathbf{E}_k$ in each element and that of $\hat{\varphi}_f$ on each face are expressed as
\begin{align}
\nonumber {\varphi _k} \simeq \sum \limits_{i = 1}^{{N_p}} {\varphi _k^i{\ell _i}}, \quad
\mathbf{E}_k  \simeq \mathop \sum \limits_{i = 1}^{{N_p}} \mathbf{E}_k^{i}{\ell _i}, \quad
\nonumber {\hat{\varphi} _f} \simeq \sum \limits_{j = 1}^{{N_{fp}}} {\hat{\varphi} _f^j{\hat{\ell} _j}}
\end{align}
where $\ell _i$, $i = 1, \ldots ,{N_p}$ , and $\hat{\ell} _j$, $j = 1, \ldots ,{N_{fp}}$, are Lagrange polynomials, $N_p$ and $N_{fp}$ are the number of nodes per element and per face, respectively, and $\varphi _k^i$, $\mathbf{E}_k^{i}$, and $\hat{\varphi} _f^j$ are the nodal values.

{\color{black}
	
Galerkin testing \eqref{weakfp0}-\eqref{weakfp1} yields the following matrix system for the local problem
\begin{align}
\label{LSlocal} & \left[ {\begin{array}{*{20}{c}}
	{\mathbb{A}_k^{\varphi \varphi}} & {\mathbb{A}_k^{\varphi \mathbf{E}}} \\ [4pt]
	{\mathbb{A}_k^{\mathbf{E} \varphi}} & {\mathbb{A}_k^{\mathbf{E} \mathbf{E}}} \\ [2pt]
	\end{array}} \right]
\left[ {\begin{array}{*{20}{c}}
	{\underline{\varphi}_k} \\ [4pt]
	{\underline{\mathbf{E}}_k} \\ [2pt]
	\end{array}} \right]
+
\left[ {\begin{array}{*{20}{c}}
	{\mathbb{A}_{k}^{\varphi \hat{\varphi}}} & \mathbb{O} \\ [4pt]
	{\mathbb{A}_{k}^{\mathbf{E} \hat{\varphi}}} & {\mathbb{A}_{k}^{\mathbf{E} \hat{\varphi}^c}} \\ [2pt]
	\end{array}} \right]
\left[ {\begin{array}{*{20}{c}}
	{\underline{\hat \varphi}_k} \\ [4pt]
	{{\hat \varphi}^c} \\ [2pt]
	\end{array}} \right]
=
\left[ {\begin{array}{*{20}{c}}
	{F_k^\varphi} \\ [4pt]
	{F_k^{\mathbf{E}}} \\ [2pt]
	\end{array}} \right]
\end{align}
where the local unknown vectors are defined as
\begin{equation}
\label{unklocal}
\underline{\varphi}_k =
\left[ {\begin{array}{*{20}{c}}
	\begin{gathered}
	{\varphi_k^1} \hfill \\
	\vdots \hfill \\
	{\varphi_k^{N_p}} \hfill \\
	\end{gathered}  
	\end{array}} \right], \quad
\underline{\mathbf{E}}_k =
\left[ {\begin{array}{*{20}{c}}
	\begin{gathered}
	{\mathbf{E}_k^1} \hfill \\
	\vdots \hfill \\
	{\mathbf{E}_k^{N_p}} \hfill \\
	\end{gathered}  
	\end{array}} \right]
\end{equation}
and $\underline{\hat{\varphi}}_k$ is a vector of dimension $N_{fe}N_{fp}\times1$, which stores the global unknowns (see \eqref{unkglobal} below) and is defined as
\begin{equation}
\label{unkmap}
\underline{\hat{\varphi}}_k = 
\left[ {\begin{array}{*{20}{c}}
	\begin{gathered}
	{\underline{\hat{\varphi}}_{k,1}} \hfill \\
	\vdots \hfill \\
	{\underline{\hat{\varphi}}_{k,{N_{fe}}}} \hfill \\
	\end{gathered}  
	\end{array}} \right], \quad
{\underline{\hat{\varphi}}_{k,{f_l}}} =
\left[ {\begin{array}{*{20}{c}}
	\begin{gathered}
	{\hat{\varphi}_{k,f_l}^1} \hfill \\
	\vdots \hfill \\
	{\hat{\varphi}_{k,f_l}^{N_{fp}}} \hfill \\
	\end{gathered}  
	\end{array}} \right].
\end{equation}
In \eqref{LSlocal}, the right hand side vectors ${F_k^{\alpha}}$, $\alpha \in \{ \varphi, \mathbf{E} \}$, correspond to the right hand sides of \eqref{weakfp0}-\eqref{weakfp1}, respectively. The matrices $\mathbb{A}_k^{\alpha \beta}$, ($\alpha \in \{ {\varphi, \mathbf{E}} \}$, $\beta \in \{ {\varphi, \mathbf{E}, \hat{\varphi}} \}$), correspond to the inner products in \eqref{weakfp0}-\eqref{weakfp1}  are the standard DG matrices, e.g., mass, stiffness, and lift matrices. For details, readers are referred to the authors' previous work~\cite{Chen2020float}. Note that $\mathbb{A}_k^{\alpha \beta}$ has dimensions $N_{\alpha} \times N_{\beta}$, where $N_{\beta}$ is the dimension of the input vector $\beta_k$ and $N_{\alpha}$ is the dimension of the output vector $\alpha_k$. 

Similarly, Galerkin testing \eqref{weakfp2}-\eqref{weakfp3} yields the following matrix system for the global problem
\begin{align}
\label{LSglobal} & \sum_{k=1}^K
\left\{
\left[ {\begin{array}{*{20}{c}}
	{\mathbb{A}_k^{\hat{\varphi} \varphi}} & {\mathbb{A}_k^{\hat{\varphi} \mathbf{E}}} \\ [4pt]
	\mathbb{O} & {\mathbb{A}_{k}^{\hat{\varphi}^c \mathbf{E}}} \\ [2pt]
	\end{array}} \right]
\left[ {\begin{array}{*{20}{c}}
	{\underline{\varphi}_k} \\ [4pt]
	{\underline{\mathbf{E}}_k} \\ [2pt]
	\end{array}} \right]
+
\left[ {\begin{array}{*{20}{c}}
	{\mathbb{A}_k^{\hat{\varphi} \hat{\varphi}}} & {\mathbb{O}} \\ [4pt]
	{\mathbb{O}} & {\mathbb{O}}\\ [2pt]
	\end{array}} \right]
\left[ {\begin{array}{*{20}{c}}
	{\underline{\hat \varphi}_k} \\ [4pt]
	{{\hat \varphi}^c} \\ [2pt]
	\end{array}} \right]
=
\left[ {\begin{array}{*{20}{c}}
	{F_k^{\hat \varphi}} \\ [4pt]
	Q^c \\ [2pt]
	\end{array}} \right] \right\}
\end{align}
where the right hand side vector ${F_k^{\hat{\varphi}}}$ corresponds to the right hand side of \eqref{weakfp2}, and the matrices $\mathbb{A}_k^{\hat{\varphi} \alpha}$, ($\alpha \in \{ {\varphi, \mathbf{E}, \hat{\varphi}} \}$), correspond to the inner products in \eqref{weakfp2}.

Solving $[\underline{\varphi}_k, \underline{\mathbf{E}}_k]^T$ in terms of $\underline{\hat{\varphi}}$ and $\hat{\varphi}^c$ from \eqref{LSlocal} yields
\begin{align}
\label{LSlocal_}
\left[ {\begin{array}{*{20}{c}}
	{\underline{\varphi}_k} \\ [4pt]
	{\underline{\mathbf{E}}_k} \\ [2pt]
	\end{array}} \right]
=
\mathbb{A}_k^{-1}
\left[ {\begin{array}{*{20}{c}}
	{F_k^\varphi} \\ [4pt]
	{F_k^{\mathbf{E}}} \\ [2pt]
	\end{array}} \right]
-
\mathbb{A}_k^{-1}
\bar{\mathbb{A}}_k
\left[ {\begin{array}{*{20}{c}}
	{\underline{\hat \varphi}_k} \\ [4pt]
	{{\hat \varphi}^c} \\ [2pt]
	\end{array}} \right]
\end{align}
where
\begin{align}
\nonumber &
\mathbb{A}_k =
\left[ {\begin{array}{*{20}{c}}
	{\mathbb{A}_k^{\varphi \varphi}} & {\mathbb{A}_k^{\varphi \mathbf{E}}} \\ [4pt]
	{\mathbb{A}_k^{\mathbf{E} \varphi}} & {\mathbb{A}_k^{\mathbf{E} \mathbf{E}}} \\ [2pt]
	\end{array}} \right], \quad
\bar{\mathbb{A}}_k =
\left[ {\begin{array}{*{20}{c}}
	{\mathbb{A}_{k}^{\varphi \hat{\varphi}}} & \mathbb{O} \\ [4pt]
	{\mathbb{A}_{k}^{\mathbf{E} \hat{\varphi}}} & {\mathbb{A}_{k}^{\mathbf{E} \hat{\varphi}^c}} \\ [2pt]
	\end{array}} \right].
\end{align}
Inserting \eqref{LSlocal_} into \eqref{LSglobal} yields a global system involving only the global unknowns
\begin{align}
\label{LSglobal_}
& \mathbb{A}_{global}
\left[ {\begin{array}{*{20}{c}}
	{\underline{\hat \varphi}} \\ [4pt]
	{{\hat \varphi}^c} \\ [2pt]
	\end{array}} \right]
=
\left[ {\begin{array}{*{20}{c}}
	{F^{\hat \varphi}} \\ [4pt]
	Q^c \\ [2pt]
	\end{array}} \right]
-
\sum_{k=1}^K {
\tilde{\mathbb{A}}_k \mathbb{A}_k^{-1}
\left[ {\begin{array}{*{20}{c}}
	{F_k^\varphi} \\ [4pt]
	{F_k^{\mathbf{E}}} \\ [2pt]
	\end{array}} \right]
}
\end{align}
where the global unknown vector is defined as
\begin{equation}
\label{unkglobal}
\underline{\hat{\varphi}} = 
\left[ {\begin{array}{*{20}{c}}
	\begin{gathered}
	{\underline{\hat{\varphi}}_1} \hfill \\
	\vdots \hfill \\
	{\underline{\hat{\varphi}}_{N_f}} \hfill \\
	\end{gathered}  
	\end{array}} \right], \quad
{\underline{\hat{\varphi}}_f} = 
\left[ {\begin{array}{*{20}{c}}
	\begin{gathered}
	{\hat{\varphi}_f^1} \hfill \\
	\vdots \hfill \\
	{\hat{\varphi}_f^{N_{fp}}} \hfill \\
	\end{gathered}  
	\end{array}} \right].
\end{equation}
and 
\begin{align}
\label{Aglobal_}
\mathbb{A}_{global} = \sum_{k=1}^K \left\{
\hat{\mathbb{A}}_k - 
{ \tilde{\mathbb{A}}_k \mathbb{A}_k^{-1} \bar{\mathbb{A}}_k }
\right\}, \quad
\tilde{\mathbb{A}}_k = 
\left[ {\begin{array}{*{20}{c}}
	{\mathbb{A}_k^{\hat{\varphi} \varphi}} & {\mathbb{A}_k^{\hat{\varphi} \mathbf{E}}} \\ [4pt]
	\mathbb{O} & {\mathbb{A}_{k}^{\hat{\varphi}^c \mathbf{E}}} \\ [2pt]
	\end{array}} \right] = \bar{\mathbb{A}}_k^T, \quad
\hat{\mathbb{A}}_k = 
\left[ {\begin{array}{*{20}{c}}
	{\mathbb{A}_k^{\hat{\varphi} \hat{\varphi}}} & {\mathbb{O}} \\ [4pt]
	{\mathbb{O}} & {\mathbb{O}}\\ [2pt]
	\end{array}} \right].
\end{align}

In~\eqref{unkmap} and~\eqref{unkglobal}, $f_l \in \{1,...,N_{fe}\}$, $f \in \{1,...,N_f\}$, $\underline{\hat{\varphi}}_{k,{f_l}}$ contains the unknowns on local face $f_l$ of element $k$, and $\underline{\hat{\varphi}}_f$ contains the unknowns on face $f$ of $\Gamma$. Apparently, each local face $f_l$ of element $k$ can be mapped to a global face $f$ of $\Gamma$. Fig.~\ref{mesh} illustrates the mapping between the nodes of the local elements (blue dots) and the nodes of the skeleton (red circles). This mapping is included in~\eqref{Aglobal_} in the summation over $k$, i.e., each local face $f_l$ of the $N_{fe}$ faces of element $k$ is mapped to one face $f$ of $\Gamma$ and, the matrix entries of the two local faces corresponding to the same $f$ are combined. The assembled matrix system from~\eqref{Aglobal_} approximately has dimensions $(N_f N_{fp}+1) \times (N_f N_{fp}+1)$. The actual size is smaller than $(N_f N_{fp}+1)$ since the nodes on $\partial \Omega$ are not included in the global problem [see \eqref{Globalfp0}]. The same mapping is done in the summation on the right hand side of~\eqref{LSglobal_}. Note that the elemental matrix ${\mathbb{A}_k^{\hat{\varphi} \varphi}}$ has dimension $N_{fe} N_{fp} \times N_p$ and the resulting vector for each $k$ has the same dimension as $\underline{\hat{\varphi}}_{k}$.
}

The size of the global system \eqref{LSglobal_} [$\sim (N_f N_{fp}+1)$] is much smaller than that of the DG method ($\sim K N_p$, see~\cite{Chen2020float}). Once $\underline{\hat{\varphi}}$ and $\hat{\varphi}^c$ are solved from \eqref{LSglobal_}, they can be used to solve $[\underline{\varphi}_k, \underline{\mathbf{E}}_k]^T$ in the local system \eqref{LSlocal_}. Since the local problems of different elements are independent from each other, they can be solved in parallel. As the dimension of \eqref{LSlocal_} is only $\sim N_p$, the computational cost of this step is relatively low and can be ignored, especially in large scale problems~\cite{Cockburn2016static}.

\section{Numerical Examples}
\subsection{Coaxial Capacitor with FPC}
\begin{figure}[!ht]
	\centering
	\subfloat[\label{Capacitora}]{\includegraphics[height=0.32\columnwidth]{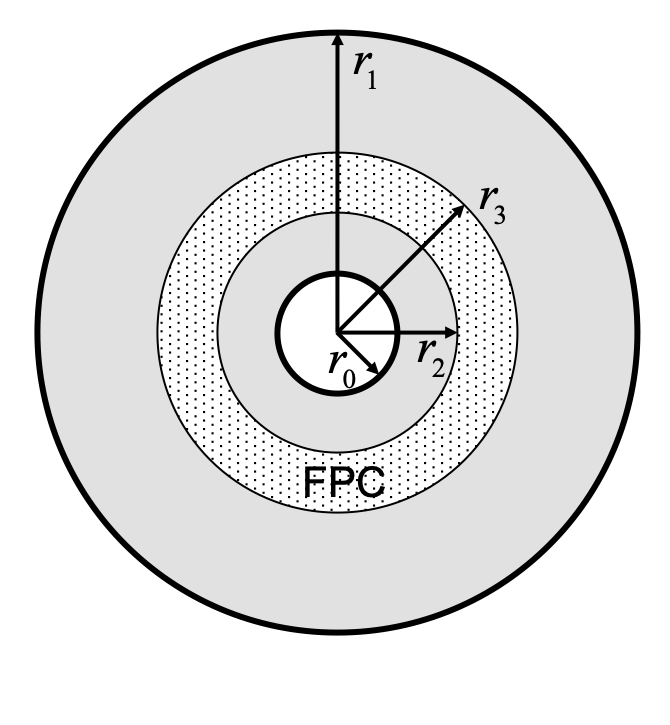}} \hspace{0.5cm}
	\subfloat[\label{Capacitorb}]{\includegraphics[height=0.32\columnwidth]{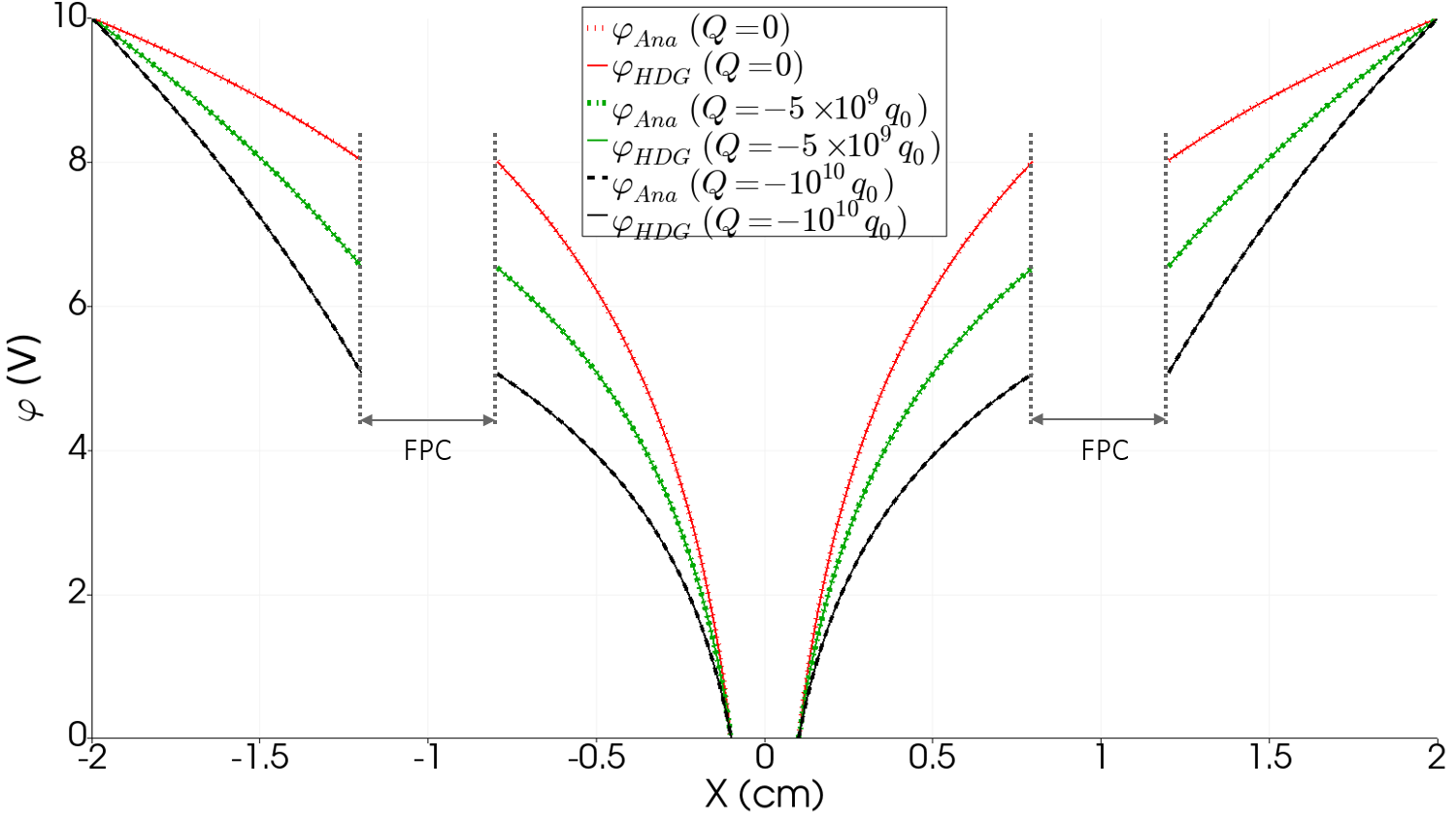}} \\
	\subfloat[\label{Capacitorc}]{\includegraphics[height=0.36\columnwidth]{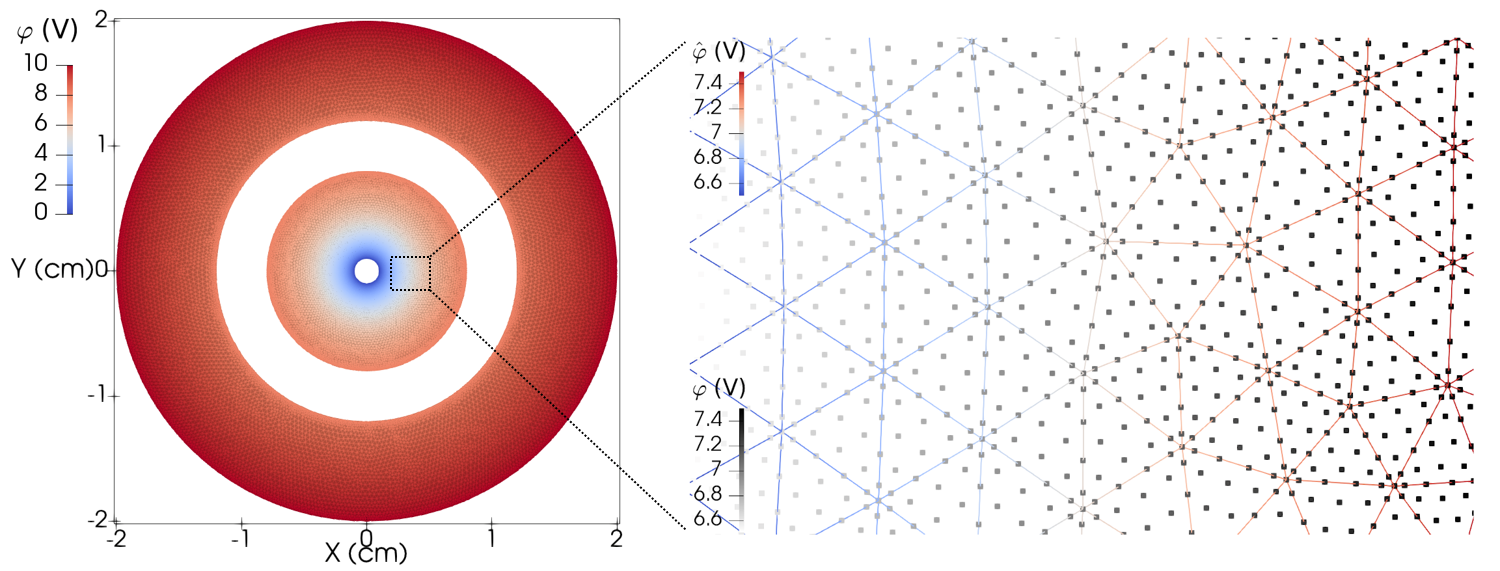}}
	\caption{(a) Schematic description of the coaxial capacitor model. (b) $\varphi$ computed by HDG and $\varphi_{Ana}$ on line $(x,y = 0)$ for different values of $Q$. (c) Illustration of the nodes where $\varphi$ and $\hat{\varphi}$ are defined.}
	\label{Capacitor}
\end{figure}

The proposed method is first validated using a canonical problem with an analytical solution. The simulation domain is illustrated in Figure~\ref{Capacitor} (a). A thin metal tube is inserted into a coaxial capacitor. The voltages applied on the inner and outer boundaries of the capacitor are $\varphi (|\mathbf{r}| = {r_0}) = {V_0}$ and $\varphi (|\mathbf{r}| = {r_1}) = {V_1}$, respectively.  The metal tube is modeled as an FPC and the FPBC is applied on $|\mathbf{r}|={r_2}$ and $|\mathbf{r}|={r_3}$. The total charge on the FPC is $Q$.The analytical solution of the electric potential is given by 
\begin{equation*}
\varphi_{Ana} (r) = \left\{ \begin{gathered}
{a_0} + {b_0}\ln (r),\;r \in [{r_0},{r_2}] \hfill \\
{a_1} + {b_1}\ln (r),\;r \in [{r_3},{r_1}] \hfill \\ 
\end{gathered}  \right.
\end{equation*}
where ${a_0} = {V_0} - {b_0}\ln ({r_0})$, ${a_1} = {V_1} - {b_e}\ln ({r_1})$, ${b_0} = {b_1} + Q/(2\pi \varepsilon )$, ${b_1} = [{V_0} - {V_0} - {C_{20}}Q/(2\pi \varepsilon )]/({C_{20}} - {C_{31}})$, and ${C_{ij}} = \ln ({r_i}/{r_j})$. In the following, ${V_0} = 0$, ${V_1} = 10$ V, ${r_0} = 0.1$ cm, ${r_1} = 2$ cm, ${r_2} = 0.8$ cm, and ${r_3} = 1.2$ cm. 

Figure~\ref{Capacitor} (b) compares the electric potential computed by HDG with $p=2$ to the analytical solution along the line $(x,y = 0)$ for $Q\in\{0, -5\times10^{10}e, -10^{10}e\}$, where $e$ is the electron charge. One can see that the numerical solution agree very well with the analytical one. The absolute value of the difference between the FPC potentials computed using HDG and the analytical solution is $1.58\times10^{-7}V$, $2.30\times10^{-8}V$ and $1.45\times10^{-8}V$ for $Q=0$, $Q=-5\times10^{9} e$ and $Q=-10^{10} e$, respectively.

\begin{table}[!ht]
\scriptsize
	\centering
	\begin{threeparttable}
		\renewcommand{\arraystretch}{1.5}
		\centering
		\caption{Dimension and condition number of the DG and HDG matrices, (wall) time and (peak) memory required by DG and HDG, and absolute error in FPC potential computed using DG and HDG for the coaxial capacitor example with zero total charge on the FPC\tnote{*}.}
		\label{nunk}
		\setlength{\tabcolsep}{3pt}
		\begin{tabular}{ p{52pt} | p{36pt} | p{36pt} | p{36pt} | p{36pt} | p{36pt} | p{36pt} | p{36pt} | p{36pt} | p{36pt} | p{36pt} }
			\hline 
			& \multicolumn{2}{c|}{$p=1$} & \multicolumn{2}{c|}{$p=2$} & \multicolumn{2}{c|}{$p=3$} & \multicolumn{2}{c|}{$p=4$} & \multicolumn{2}{c}{$p=5$} \\ \hline
			& DG & HDG & DG & HDG & DG & HDG & DG & HDG & DG & HDG \\ \hline
			Dimension   & 254,838 & 252,319 & 509,676 & 378,478 & 849,460 & 504,637 & 1,274,190 & 630,796 & 1,783,866 & 756,955 \\ \hline
			Condition \# & 1.33$\times 10^8$ & 1.17$\times 10^8$ & 5.44$\times 10^8$ & 1.67$\times 10^8$ & 16.17$\times 10^8$ & 2.62$\times 10^8$ & 39.4$\times 10^8$ & 3.47$\times 10^8$ & 84.2$\times 10^8$ & 4.57$\times 10^8$ \\ \hline			
			Time (s) & 1.97 & 1.79 & 5.35 & 3.83 & 10.5 & 7.09 & 18.8 & 9.29 & 32.3 & 12.6 \\ \hline
			Memory (GB)& $0.41$ & $0.38$ & $1.07$ & $0.93$ & $2.11$ & $1.66$ & $3.89$ & $2.42$ & $6.05$ & $3.21$ \\ \hline
			$\mathrm{Error}$ (V) & $2.86$$\times10^{-4}$ & $2.82$$\times10^{-4}$ & $2.30$$\times10^{-7}$ & $2.26$$\times10^{-7}$ & $2.01$$\times10^{-7}$ & $1.99$$\times10^{-7}$ & $1.90$$\times10^{-7}$ & $1.85$$\times10^{-7}$ & $1.90$$\times10^{-7}$ & $1.83$$\times10^{-7}$ \\ \hline
		\end{tabular}
		\smallskip
		\scriptsize
		\begin{tablenotes}
			\item[*] {The matrix systems are solved using UMFPACK (multifrontal sparse LU factorization) implemented by Matlab\textsuperscript{\textregistered} on a workstation with Intel\textsuperscript{\textregistered} Xeon\textsuperscript{\textregistered} E5-2680 v4 processor (2.40 GHz base clock, 35 MB cache, and 14 cores/28 threads). Only 24 threads are used.}
		\end{tablenotes}
	\end{threeparttable}
\end{table}
Figure~\ref{Capacitor} (c) shows the nodes where $\varphi$ and $\hat{\varphi}$ are defined for $p=5$. The degrees of freedom of $\hat{\varphi}$ only correspond to the nodes on the wireframe while those of $\varphi$ correspond to all nodes. Same as DG, the nodes of $\varphi$ are doubly defined on the wireframe. 

Table~\ref{nunk} presents the dimension and the condition number of the DG~\cite{Chen2020float} and HDG matrices, the (wall) time and the (peak) memory required by DG and HDG, and the absolute error in FPC potential computed using DG and HDG. The table clearly shows that the difference in matrix dimensions increase with increasing $p$. This is easy to see from Figure~\ref{Capacitor} (c). With increasing $p$, the ratio of the number of nodes in the interior of each element to that on the element surface increases. Accordingly, as also shown by the table, the benefits introduced by HDG in wall time and peak memory requirement become more significant with increasing $p$.
{\color{black}
Table~\ref{nunk} also shows that, as expected, the condition numbers of the DG and HDG matrices increase with increasing $p$~\cite{Hesthaven2008}. However, the condition number of the HDG matrix grows much slower. This is because the dimension of the mesh skeleton is always one smaller than that of the mesh itself. Finally, Table~\ref{nunk} shows that the absolute error in FPC potential computed using DG and HDG is almost the same. It should be noted here that, in this example, the accuracy is limited by the representation of a curved structure using linear elements. Hence, the error does not decrease much as $p$ is increased beyond $2$.

It is worth to include here a discussion comparing the proposed HDG method to the conventional FEM-based methods, namely MRM~\cite{Dong2010} and VPM~\cite{Konrad1996}. As briefly discussed in Section~\ref{sec1}, MRM and VPM are applicable only when $Q=0$. For this case, it has been shown in~\cite{Chen2020float} that DG has the same accuracy as MRM and is more accurate than VPM. But the dimension of the DG matrix system is larger than those of MRM and VPM due to the duplication of the nodes on the element surfaces. This difference becomes smaller for larger values of $p$, for which the number of the interior nodes is larger than the number of the nodes on the element surface. The same observation regarding the accuracy and the matrix dimension holds true when the proposed HDG method is compared to MRM and VPM. As shown in Table~\ref{nunk}, the reduction in the matrix dimension obtained by using HDG becomes more significant with increasing $p$. However, in general, HDG still needs to solve a larger matrix system than MRM and VPM. The efficiency of HDG becomes comparable to that of FEM only when higher-order basis functions are used~\cite{Kirby2012ToCGorHDG2D, Yakovlev2016ToCGorHDG3D}. Detailed studies in~\cite{Kirby2012ToCGorHDG2D, Yakovlev2016ToCGorHDG3D} show that the efficiency also strongly depends on the specific problem being analyzed and the mesh used for discretizing the problem since the conditioning of the matrix system plays an important role. The efficiency of HDG can be improved using h/p- refinement techniques and employing a non-conformal mesh when possible. Note that the implementation of these approaches is significantly easier for HDG (and DG) than FEM.
}

\subsection{Plasmonic-enhanced Photoconductive Antenna}
\begin{figure}[!bp]
	\centering
	\subfloat[\label{PCAa}]{\includegraphics[height=0.38\columnwidth]{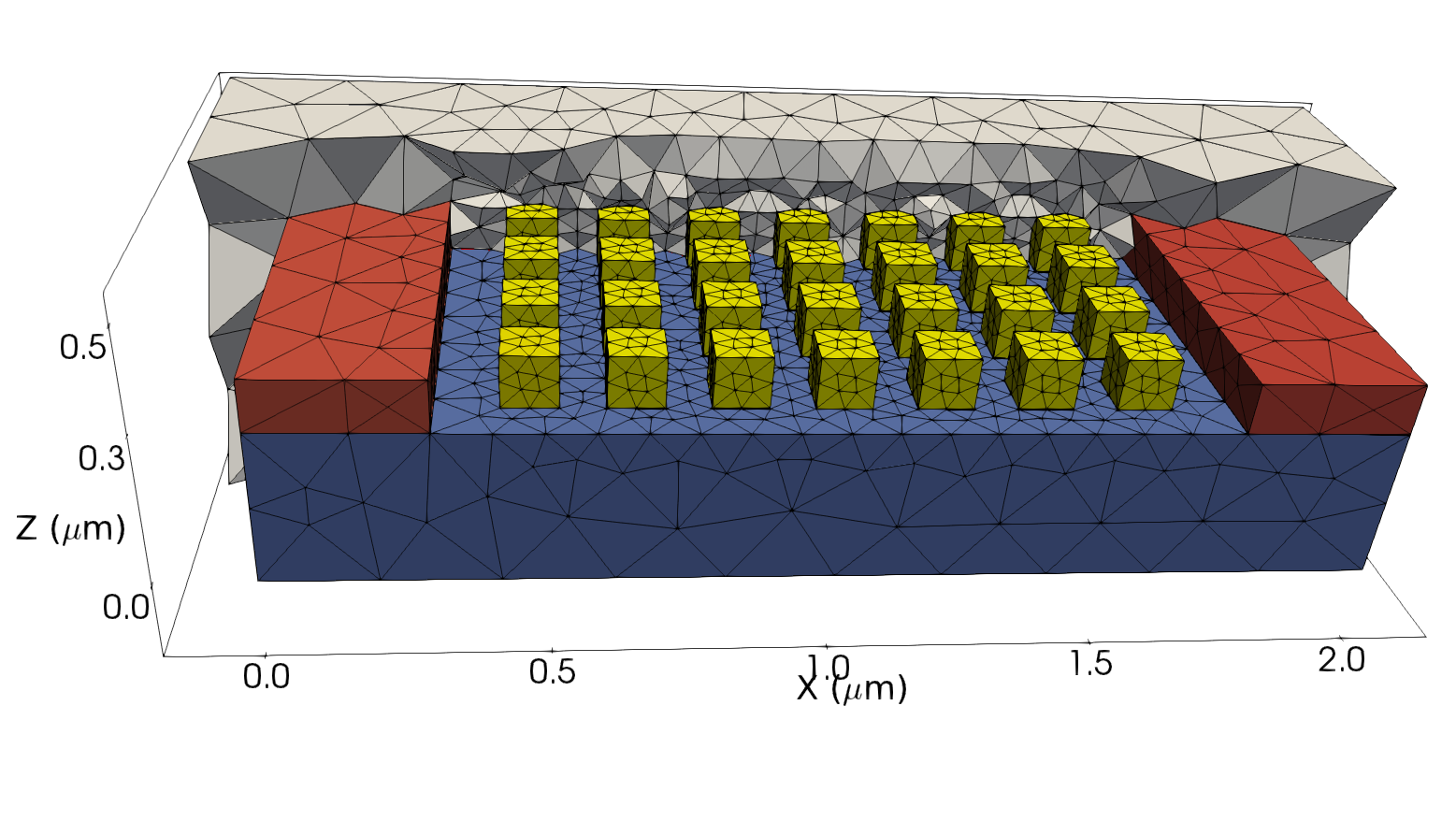}} \\
	\subfloat[\label{PCAv}]{\includegraphics[height=0.4\columnwidth]{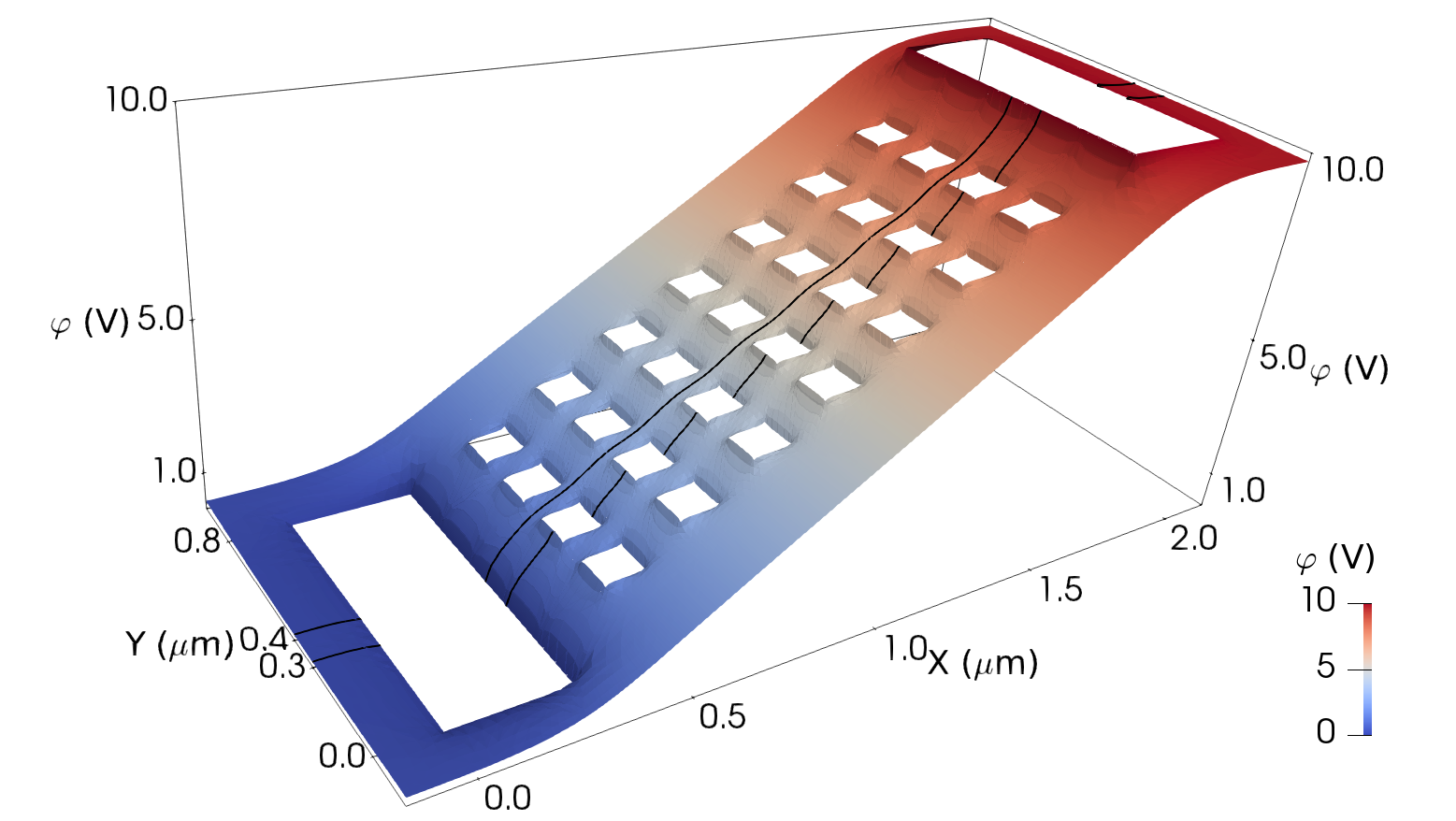}} \\
	\subfloat[\label{PCAc}]{\includegraphics[height=0.3\columnwidth]{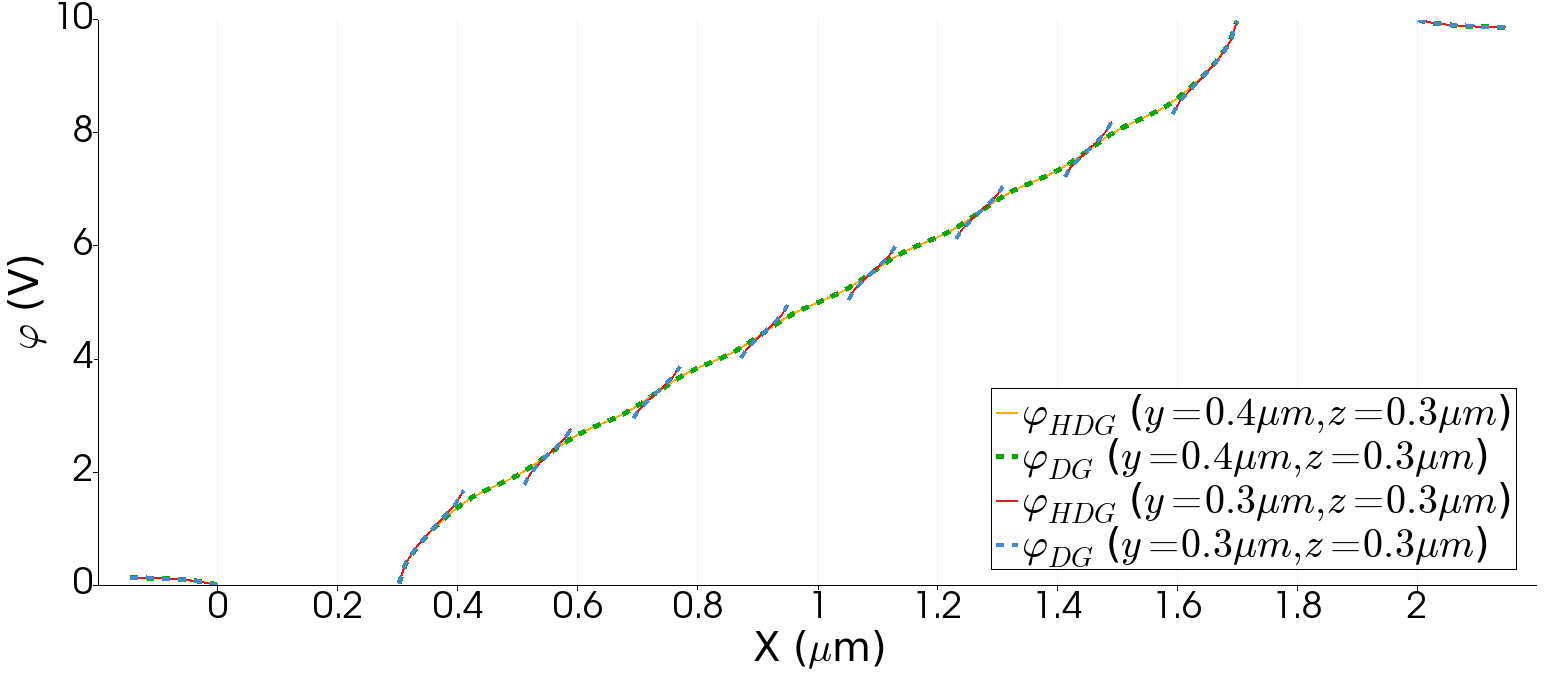}}
	\caption{(a) Geometry of the plasmonic-enhanced PCA. (b) $\varphi$ computed using HDG on plane $(x,y,z=0.3~\mu {\rm m})$. (c) $\varphi$ computed by HDG and DG on lines $(x, y=0.3 ~\mu\mathrm{m}, z = 0.3~\mu\mathrm{m})$ and $(x, y=0.4 ~\mu\mathrm{m}, z = 0.3~\mu\mathrm{m})$. The positions of these lines are marked by the black lines in (b). }
	\label{PCA}
\end{figure}

Next, a plasmonic-enhanced photoconductive antenna (PCA) is considered. The device geometry is shown in Fig.~\ref{PCA} (a). The semiconductor layer (blue) is made of GaAs, which is a photoconductive material that can absorb optical electromagnetic (EM) wave energy and generate terahertz (THz) signals~\cite{Chen2019multiphysics}. The metallic nanostructures (yellow) are designed to enhance the local EM fields and hence increase the optical-to-THz efficiency~\cite{Park2012}. Throughout the operation of this device, a bias voltage is applied on the two electrodes (red), generating a static electric field. Thus, to model this device, one needs to solve the electrostatic problem under the bias voltage, aside from the transient EM response~\cite{Chen2020steadystate, Chen2020float, Chen2019multiphysics}.

The semiconductor layer has relative permittivity of $10.9$ and the surrounding area (gray) is air. The computation domain is truncated at the outmost boundaries with homogeneous Neumann boundary condition. Dirichlet boundary conditions are applied on the surfaces of the electrodes, with $\varphi(\mathbf r)=0$ on the left one (cathode) and $\varphi(\mathbf r)=10 $~V on the right one (anode). Since the metallic nanostructures are isolated conductors, they act as FPCs and an independent FPBC is applied on each block of the nanostructures.

Fig.~\ref{PCA} (b) shows $\varphi$ and $\hat{\varphi}$ on the plane $(x,y,z=3 ~\mu \mathrm{m})$ solved from the proposed HDG method with $p=4$. As expected, $\varphi(\mathbf r)$ is constant on the surface of each FPC block. Meanwhile, the potential values on different FPC blocks are different because the FPCs are isolated from each other. One consequence of the inhomogeneous electric potential distribution is that strong local static electric fields are generated near the FPC, which greatly influences the carrier mobilities in the semiconductor layer~\cite{Chen2020steadystate} and hence influences the device performance~\cite{Moon2012}.

\begin{table}[!b]
\scriptsize
	\centering
	\begin{threeparttable}
		\renewcommand{\arraystretch}{1.5}
		\centering
		\caption{Dimension and condition number of the DG and HDG matrices, and (wall) time and (peak) memory required by DG and HDG for the plasmonic-enhanced PCA example\tnote{*}.}
		\label{nunk1}
		\setlength{\tabcolsep}{3pt}
		\begin{tabular}{ p{52pt} | p{36pt} | p{36pt} | p{36pt} | p{36pt} | p{36pt} | p{36pt} | p{36pt} | p{36pt} }
			\hline 
			& \multicolumn{2}{c|}{$p=3$} & \multicolumn{2}{c|}{$p=4$} & \multicolumn{2}{c|}{$p=5$} & \multicolumn{2}{c}{$p=6$} \\ \hline
			& DG & HDG & DG & HDG & DG & HDG & DG & HDG \\ \hline
			Dimension   & 81,440 & 75,908 & 142,520 & 113,848 & 228,032 & 159,376 & 342,048 & 212,492 \\ \hline	
			Condition \# 
			& 3.39$\times10^8$      & 3.16$\times10^6$ 
			& 5.08$\times10^9$    & 4.95$\times10^6$
			& 2.37$\times10^{10}$   & 6.68$\times10^6$
			& 1.87$\times10^{11}$  & 8.36$\times10^6$ \\ \hline
			Time (s) & 1.57   & 1.46  & 3.38  & 2.64  & 7.80  & 5.88  & 16.3 & 10.3 \\ \hline
			Memory (GB)   & $0.91$ & $0.88$ & $2.24$ & $1.90$ & $4.91$& $3.67$& $9.63$ & $6.47$  \\ \hline
		\end{tabular}
		\smallskip
		\scriptsize
		\begin{tablenotes}
			\item[*] {The matrix systems are solved using Intel\textsuperscript{\textregistered} MKL PARDISO (parallel direct sparse solver) (v2018.2) on a workstation with Intel\textsuperscript{\textregistered} Xeon\textsuperscript{\textregistered} E5-2680 v4 processor (2.40 GHz base clock, 35 MB cache, and 14 cores/28 threads). Only 24 threads are used. }
		\end{tablenotes}
	\end{threeparttable}
\end{table}
Fig.~\ref{PCA} (c) shows the solutions on the lines $(x, y=0.3 ~\mu\mathrm{m}, z = 0.3~\mu\mathrm{m})$ and $(x, y=0.4 ~\mu\mathrm{m}, z = 0.3~\mu\mathrm{m})$ computed by HDG and DG~\cite{Chen2020float} with $p=4$. The maximum difference between the solutions obtained using the two solvers is $1.1\times 10^{-4}~\mathrm{V}$. Note that no volumetric meshes are used inside the FPCs and the electrodes since they are treated as boundary conditions by HDG and DG. For practical device simulations, this treatment can save considerable amount of computational resources since finer meshes are usually required near the nanostructures~\cite{Chen2020steadystate, Chen2019multiphysics}.

{\color{black}
Table~\ref{nunk1} presents the dimension and the condition number of the DG and HDG matrices, and the (wall) time and the (peak) memory required by DG and HDG. Just like the previous example, the reduction in the matrix dimension and accordingly savings in simulation time and memory requirement, which are obtained by using HDG, become more significant with increasing $p$. Both condition numbers increase with increasing $p$, both the growth rate of the condition number of the HDG matrix is much smaller. It should be noted here that, compared to the previous example, the higher condition numbers (for both DG and HDG) are due to the highly skewed mesh elements. This could be avoided by using small elements throughout the whole computation domain but this comes with increased simulation time and memory requirement. }

\begin{figure}[!b]
	\centering
	\subfloat[\label{Arrestera}]{\includegraphics[height=0.42\columnwidth]{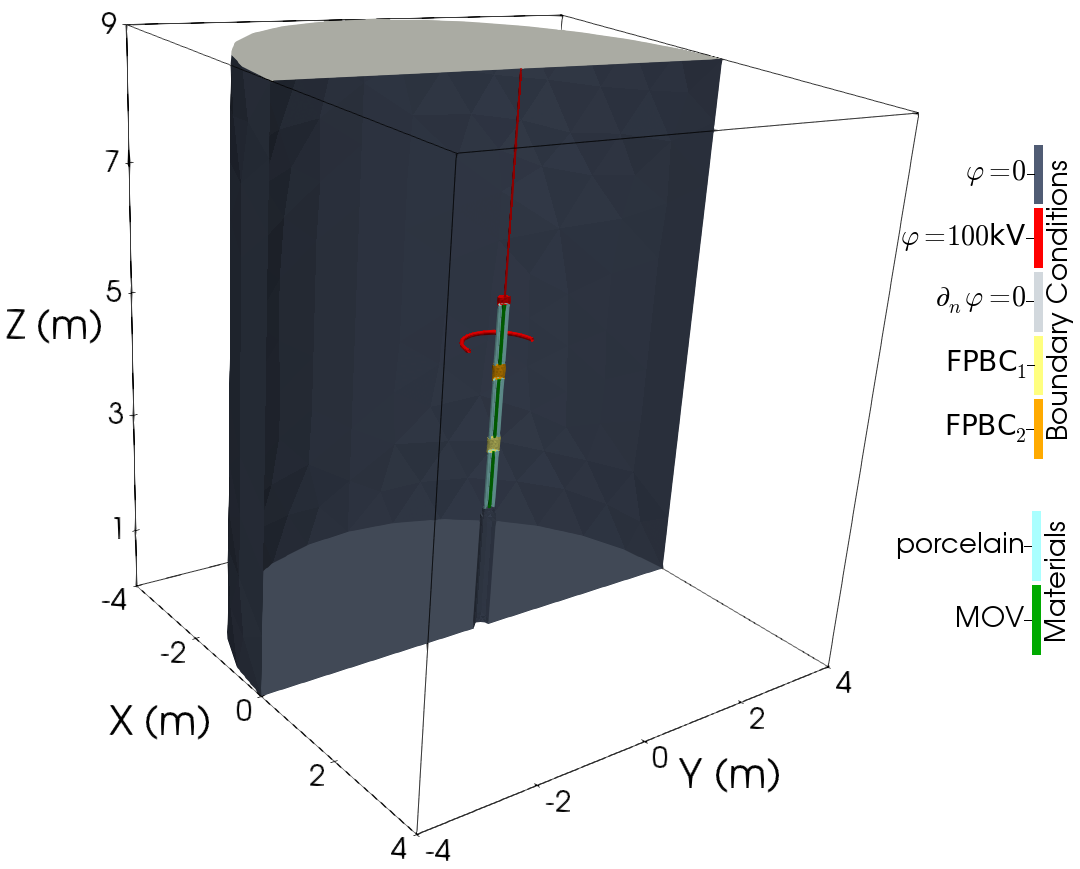}}
	\hspace{0.2cm}
	\subfloat[\label{Arresterb}]{\includegraphics[height=0.42\columnwidth]{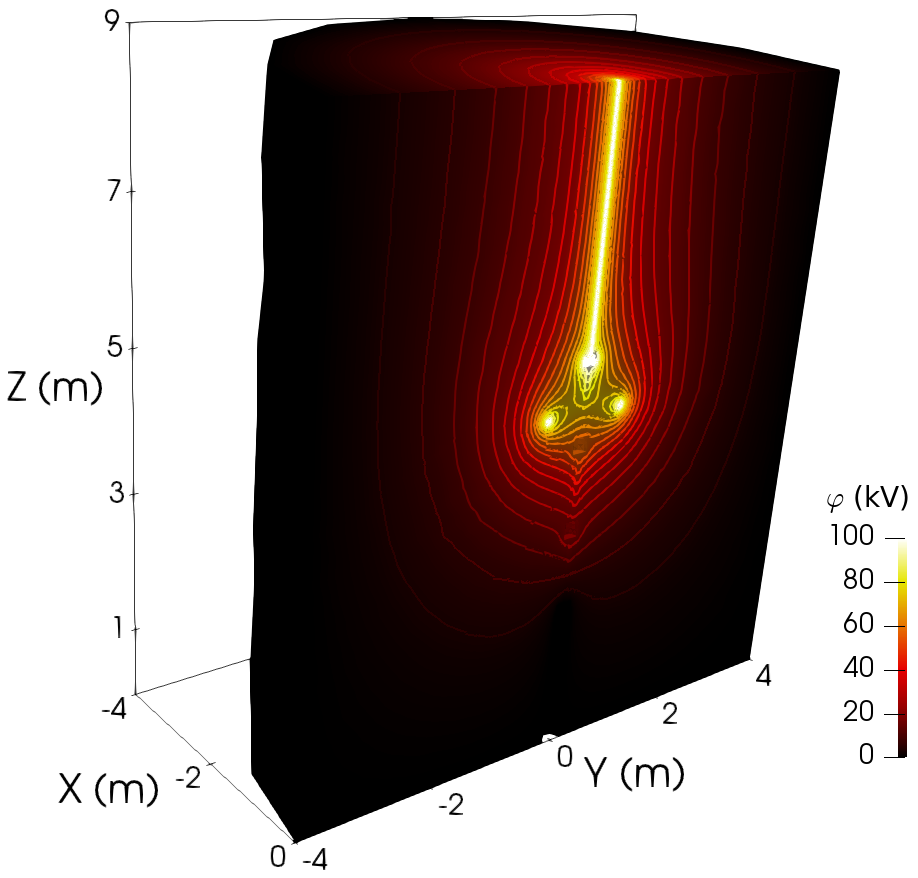}} \\
	%	\hspace{0.15cm}
	\subfloat[\label{Arresterc}]{\includegraphics[width=0.65\columnwidth]{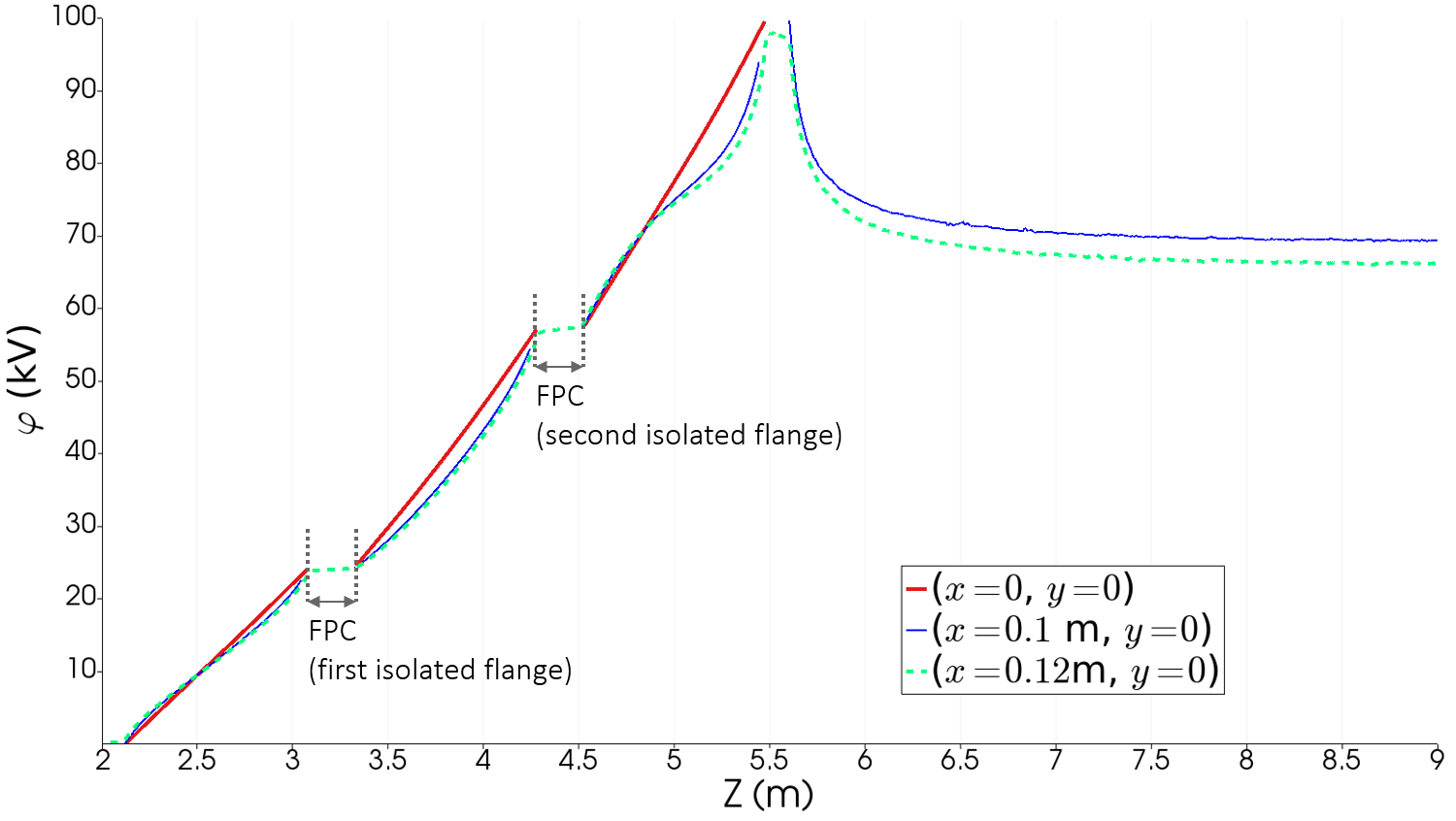}}
	\caption{(a) Geometry of the IEC surge arrester. (b) $\varphi$ computed by HDG. (c) $\varphi$ computed by HDG on lines $(x=0, y=0, z)$, $(x=0.1 ~\mathrm{m}, y=0, z)$ and $(x=0.12 ~\mathrm{m}, y=0, z)$.}
	\label{Arrester}
\end{figure} 

\subsection{Surge Arrester}
Next, the proposed method is used to compute the electric potential on an IEC surge arrester~\cite{Arrester2006}. The model is shown in Figure~\ref{Arrester} (a). The arrester consists of three segments of metal-oxide varistor (MOV) column. Each segment is surrounded by a porcelain layer. The pedestal and the surrounding cylinder are grounded ($\varphi=0$) and a high voltage ($\varphi=100\mathrm{kV}$) is applied on the lead and the grading ring~\cite{Arrester2006}. The MOV columns are separated by two isolated metal flanges, which are considered as FPCs and are modeled with two independent FPBCs in the proposed HDG method. 

The diameters of the MOV, the inner wall of the porcelain layer, the outer wall of the porcelain layer, the flanges, the pedestal, and the lead are $60~\mathrm{mm}$, $140~\mathrm{mm}$, $200~\mathrm{mm}$, $230~\mathrm{mm}$, $280~\mathrm{mm}$, and $40~\mathrm{mm}$, respectively. The major and minor diameters of the grading ring are $1130~\mathrm{mm}$ and $70~\mathrm{mm}$, respectively. The heights of the pedestal, each metal flange, each MOV segment, and the lead are $2000~\mathrm{mm}$, $120~\mathrm{mm}$, $960~\mathrm{mm}$, and $3400~\mathrm{mm}$, respectively. The diameter and height of the surrounding cylinder are $8000~\mathrm{mm}$ and $9000~\mathrm{mm}$, respectively, which are determined by the minimum phase-to-earth clearance~\cite{Arrester2006}. The relative permittivity of the MOV column and  the porcelain layer are $\epsilon_r=800$ and $\epsilon_r=5$, respectively.

Figure~\ref{Arrester} (b) shows the electric potential distribution computed with the HDG scheme using $p=4$. Figure~\ref{Arrester} (c) shows the solutions on the lines $(x=0, y=0, z)$, $(x=0.1 ~\mathrm{m}, y=0, z)$ and $(x=0.12 ~\mathrm{m}, y=0, z)$. The potential values on the two FPCs are $24.15\mathrm{kV}$ and $54.17\mathrm{kV}$. These results agree with the data reported in~\cite{Arrester2006,Amann2014,Andjelic2018}. The maximum difference between the solutions of HDG and DG over the whole domain is $1.1\times10^{-5}\mathrm{kV}$.
{\color{black}

Table~\ref{nunk2} presents the dimension and the condition number of the DG and HDG matrices, and the (wall) time and the (peak) memory required by DG and HDG. Same conclusions as the previous two examples apply for the results presented in this table.}

\begin{table}[!ht]
\scriptsize
	\centering
	\begin{threeparttable}
		\renewcommand{\arraystretch}{1.5}
		\centering
		\caption{Dimension and condition number of the DG and HDG matrices, and (wall) time and (peak) memory required by DG and HDG for the surge arrester example\tnote{*}.}
		\label{nunk2}
		\setlength{\tabcolsep}{3pt}
		\begin{tabular}{ p{52pt} | p{36pt} | p{36pt} | p{36pt} | p{36pt} | p{36pt} | p{36pt} }
			\hline 
			& \multicolumn{2}{c|}{$p=3$} & \multicolumn{2}{c|}{$p=4$} & \multicolumn{2}{c}{$p=5$} \\ \hline
			& DG & HDG & DG & HDG & DG & HDG \\ \hline
			Dimension   & 402,740 & 376,802 & 704,795 & 565,202 & 1,127,672 & 791,282 \\ \hline	
			Condition \# 
			& 9.69$\times 10^9$     & 6.97$\times 10^7$ & 6.83$\times 10^{10}$ & 1.00$\times 10^8$ & 7.26$\times 10^{11}$  & 1.48$\times 10^8$ \\ \hline
			Time (s) & 8.29 & 5.63 & 20.3 & 13.5 & 43.6 & 30.1 \\ \hline
			Memory (GB)   & $4.31$ & $4.25$ & $11.0$ & $9.46$ & $23.9$ & $18.2$ \\ \hline
		\end{tabular}
		\smallskip
		\scriptsize
		\begin{tablenotes}
			\item[*] {The matrix systems are solved using Intel\textsuperscript{\textregistered} MKL PARDISO (parallel direct sparse solver) (v2018.2) on a workstation with Intel\textsuperscript{\textregistered} Xeon\textsuperscript{\textregistered} E5-2680 v4 processor (2.40 GHz base clock, 35 MB cache, and 14 cores/28 threads). Only 24 threads are used.}
		\end{tablenotes}
	\end{threeparttable}
\end{table}

\section{Conclusions}\label{sec5}
A HDG scheme for modeling FPCs in electrostatic problems is developed. The local problem is formulated as a Dirichlet BVP, the global problem is formulated in the unknown electric potential and the unknown floating potential values of each FPC. The proposed HDG scheme retains the advantages of the DG scheme previously proposed for FPC modeling in electrostatic simulations, i.e., FPCs that can account for non-zero charge conditions, accurate solution, and ease of implementation in an existing code. Meanwhile, it significantly reduces the number of degrees of freedom as compared to DG, leading to a reduced computational cost. 

%\backmatter

\section*{Acknowledgment}
This work is supported in part by the King Abdullah University of Science and Technology (KAUST) Office of Sponsored Research (OSR) under Award No 2016-CRG5-2953 and in part by the National Natural Science Foundation of China under Grant 61701424. The authors would like to thank the King Abdullah University of Science and Technology Supercomputing Laboratory (KSL) for providing the required computational resources.

\end{document}